\documentclass[12pt]{amsart}

\def\@typesizes{%
       \or{5}{6.5}\or{6}{7.5}\or{7}{8.5}\or{8}{11}\or{9}{12}%
       \or{10}{13}
       \or{\@xipt}{14}\or{\@xiipt}{15}\or{\@xivpt}{18}%
       \or{\@xviipt}{20}\or{\@xxpt}{24}}

\usepackage{amsthm}
\usepackage{amsmath}
\usepackage{amssymb}
\usepackage{graphicx}
\usepackage{enumerate}
\usepackage{color}
\allowdisplaybreaks
\numberwithin{equation}{section}
\numberwithin{figure}{section}

\theoremstyle{plain}
\newtheorem{theorem}{ Theorem}[section]
\newtheorem{proposition}[theorem]{ Proposition}
\newtheorem{lemma}[theorem]{ Lemma}
\newtheorem{corollary}[theorem]{ Corollary}
\newtheorem{example}[theorem]{ Example}
\newtheorem{remark}[theorem]{ Remark}
\newtheorem{definition}[theorem]{ Definition}
\newtheorem{conjecture}{ Conjecture}

\def\BET{\begin{theorem}}
\def\ENT{\end{theorem}}
\def\BEP{\begin{proposition}}
\def\ENP{\end{proposition}}
\def\BEL{\begin{lemma}}
\def\ENL{\end{lemma}}
\def\BEC{\begin{corollary}}
\def\ENC{\end{corollary}}
\def\BEE{\begin{example} \rm}
\def\ENE{\end{example}}
\def\BER{\begin{remark} \rm}
\def\ENR{\end{remark}}
\def\BED{\begin{definition} \rm}
\def\END{\end{definition}}
\def\BECJ{\begin{conjecture}}
\def\ENCJ{\end{conjecture}}

\def\bea{\begin{eqnarray}}
\def\eea{\end{eqnarray}}

\def\beas{\begin{eqnarray*}}
\def\eeas{\end{eqnarray*}}

\def\beq{\begin{equation}}
\def\eeq{\end{equation}}

\def\beal{\begin{align*}}

\def\eeal{ \end{align*} }

\def\roweq{\nonumber \\ &=& }
\def\rowleq{\nonumber \\  & \leq & }

\def\bbC{{\mathbb C}}
\def\bbD{{\mathbb D}}

\def\bbN{{\mathbb N}}

\def\bbR{{\mathbb R}}

\def\cP{{\mathcal P}}

\def\cW{{\mathcal W}}

\begin{document}

\title{On solid cores and hulls of  weighted Bergman spaces $A_{\mu}^1$}

\maketitle

\noindent Jos\'e Bonet

\noindent
Instituto Universitario de Matem\'{a}tica Pura y Aplicada IUMPA,
Universitat Polit\`{e}cnica de Val\`{e}ncia,  E-46071 Valencia, Spain

\bigskip

\noindent Wolfgang Lusky

\noindent Institut f\"ur Mathematik, Universit\"at Paderborn, D-33098 Paderborn, Germany


\bigskip

\noindent  Jari Taskinen (corresponding author)

\noindent Department of Mathematics and Statistics, P.O. Box 68,
University of Helsinki, 00014 Helsinki, Finland

email: jari.taskinen@helsinki.fi


\bigskip


\smallskip

\noindent {\small \it To appear in Journal of Mathematical Sciences}

\bigskip

Abstract.
We consider weighted Bergman spaces $A_\mu^1$ on the unit disc as well as
the corresponding spaces of entire functions, defined using
non-atomic Borel measures with radial symmetry. By extending the techniques
from the case of reflexive Bergman spaces we characterize the solid  core
of $A_\mu^1$. Also, as a consequence of a characterization of solid
$A_\mu^1$-spaces we show that, in the case of entire functions, there indeed
exist solid  $A_\mu^1$-spaces. The second part of the paper
is restricted to the case of the unit disc and it contains a characterization
of the solid hull of $A_\mu^1$, when $\mu $ equals the weighted Lebesgue measure
with weight $v$. The results are based on a duality relation
of weighted  $A^1$- and $H^\infty$-spaces, the validity of which requires the
assumption that  $- \log v$ belongs to the class $\cW_0$, studied in a number of
publications; moreover, $v$ has to satisfy condition $(b)$, introduced by the
authors. The exponentially decreasing weight
$v(z) = \exp( -1 /(1-|z|)$ provides an example satisfying both assumptions.

\section{Introduction and preliminaries}
\label{sec0}

The solid hulls and cores of spaces of analytic functions on the unit disc $\bbD= \{ z \in \mathbb{C} : |z| < 1  \}$ or the entire plane
$\bbC$ have been investigated by many authors. We refer the reader to  the recent books \cite{JVA} and \cite{Pav-book} and the many references therein.
In the series of papers \cite{BT}--\cite{BLT2} the authors have presented the solid hulls and cores
of the weighted $H^\infty$-spaces $H^\infty_v$ on $\bbD$ or $\bbC$ for a large class of radial
weights $v$ as well as their
Bergman space analogues $A_\mu^p$ for $1 < p < \infty$. Earlier, the cases of standard
weights and $d \mu(r) = (1-r)^{\alpha}dr$, $\alpha > 0$, were considered in
\cite{AS} and  \cite{JVA}.

In this note we want to extend the results of \cite{BLT2} to weighted Bergman
spaces $A_{\mu}^p$ for $p=1$. The spaces are defined on the unit disc
$\bbD$ or on the entire  plane. ({\it Fock spaces} are usually considered as the
Bergman space analogues of spaces of entire functions, but these are defined with Gaussian weight functions, which is not required here.
Thus, we keep here the term Bergman space also for entire functions.)
Consider $R =1 $ or $R = \infty$. We study   holomorphic functions
$ f: R \cdot \mathbb{D} \rightarrow
\mathbb{C}$ where $R \cdot \mathbb{D} =\mathbb{D}$ if $R=1$  and $R  \cdot \mathbb{D}= \mathbb{C}$ if $R = \infty$. Let $ \hat{f}(k)$ be the Taylor coeffients of $f$, i.e. $ f(z) = \sum_{k=0}^{\infty} \hat{f}(k) z^k$. We
take a non-atomic positive bounded Borel measure $ \mu$ on $[0,R[$  such that $\mu([r, R[) > 0$ for every
$ r > 0$ and $\int_0^R r^n d \mu(r) < \infty$ for all $n> 0$. Put, for $1 \leq p < \infty$,
\[
\Vert f \Vert_p = \left(\frac{1}{2 \pi}\int_0^R \int_0^{2 \pi}|f(re^{i \varphi})|^p d \varphi d \mu(r) \right)^{1/p}
\]
and let
\[
A_{\mu}^p = \{ f: R \cdot \mathbb{D} \rightarrow \mathbb{C} : f \mbox{ holomorphic with  }  \Vert f \Vert_p < \infty \}.
\]
We will also consider the weighted spaces
\[
H_v^\infty = \{ f : \bbD \to \bbC  \mbox{ holomorphic} :
\Vert f \Vert_v = \sup\limits_{z \in \bbD} v(z) |f(z)| < \infty \} ,
\]
where the weight $v: \bbD \to (0, \infty) $ is a continuous and radial
($v(z) =v(|z|)$) function which is decreasing with respect to $r = |z|$, and
$\lim_{r \to 1^-} v(r)  = 0$.

Let $A$ be a vector space of holomorphic functions on $ R \cdot \mathbb{D}$ containing the polynomials.
The \textit{solid core} is defined as
\[
s(A) =  \{ f \in A : g \in A \mbox{ for all holomorphic } g
\mbox{ with } |\hat{g}(k)| \leq |\hat{f}(k) | \mbox{ for all } k \}
\]
and the \textit{solid hull} as
\[
S(A) = \{ g: \mathbb{D} \rightarrow \mathbb{C}  \mbox{ holomorphic :
there is } f \in A  \mbox{ with } |\hat{g}(k)| \leq |\hat{f}(k) | \mbox{ for all } k \}.
\]
The space $A$ is called \textit{solid} if $A= S(A)$. The concept of a solid
hull will also be discussed in the beginning of Section \ref{sec2}.

Here, in Theorem \ref{th1.2} we will transfer Theorem 4.1. of \cite{BLT2} to the case $p=1$.
This result concerns the characterization of solid Bergman spaces $A_{\mu}^1$, and it is motivated
by the fact that such spaces indeed exist in the case $R= \infty$ (only), as will be shown
in Example \ref{ex1} and Corollary \ref{cor1.5}. In Theorem \ref{th1.6} we determine the solid
cores for all Bergman spaces $A_{\mu}^1$.

We also present in Section \ref{sec2} how duality theory can be used for new
results on certain solid hulls; see the beginning of Section \ref{sec2} for
detailed definitions. In particular, 
we construct the solid hull $S_{BK} (A_\mu^1)$ of $A_\mu^1$ for $R=1$ by using the
known solid core of the space $H^\infty_v$ in \cite{BLT1}. This result is more
special than the above one for solid cores, since we need to restrict to the
case $\mu$  is the weighted Lebesgue measure $d \mu = v dA = v \pi^{-1} r dr d
\varphi$, where the weight $v$ needs to satisfy some special assumptions in
addition to those mentioned above. Examples
of such weights include important cases like exponentially decreasing weights.

For a holomorphic $g$ and $r > 0$ we define
\[
M_p(g,r) = \left(\frac{1}{2 \pi} \int_0^{2 \pi} |g(r e^{i \varphi})|^pd\varphi
\right)^{1/p}
\]
and denote the Dirichlet projections by $P_n g(z) = \sum_{k=0}^n \hat{g}(k) z^k$, $n \in \bbN$. It is well-known that, for $1 < p < \infty$, there are
constants $c_p>0$, not depending on $g$, $n$ or $r$, such that
$ M_p(P_ng,r) \leq c_p M_p(g,r)$. Moreover we have $\lim_{n \rightarrow \infty} M_p(g - P_ng, r) =0$. Hence we obtain
\[
\Vert P_nf \Vert_p \leq c_p  \Vert f \Vert_p \mbox{ for all } f \in   A_{\mu}^p \mbox{ and all } n \mbox{ and } \lim_{n \rightarrow \infty}  \Vert f - P_nf \Vert_p =0.
\]
In particular we see that the monomials $z \mapsto z^n$, $n \in \bbN_0 = \{ 0,1,2, \ldots
\} = \bbN \cup \{0 \}$, form a (Schauder) basis of
$A_{\mu}^p$ if $1 < p < \infty$. On the other hand, denoting by $H^1$
the Hardy space of all holomorphic functions on
$\mathbb{D}$ which are bounded under $\sup_{0 \leq r <1}M_1(\cdot, r)$,
it is well known that  the operator norm of $P_n : H^1 \to H^1$
tends to infinity  as $n \rightarrow \infty$.  Details can be seen in \cite{PD}
and \cite{W}. For the terminology and definitions on bases in Banach spaces, see also  \cite{LT}.

In the rest of the article $[r]$ denotes the largest integer smaller or equal
than $r>0$. By $c, c_1, c_2, C, C'$ etc. we denote generic positive constants,
the actual value of which may vary from place to place.

\section{Solid core and examples of solid $A_{\mu}^1$-spaces.}
\label{sec1}

In this section we extend  Theorem 4.1. of \cite{BLT2} concerning
the characterization of solid Bergman spaces to the case $p=1$ and
also determine the solid cores for all spaces $A_{\mu}^1$.
We consider both cases $R=1$ or $R= \infty$ unless otherwise specified. At
first we recall a fundamental result from \cite{ahl}, which  concerns
equivalent representations of the norm of the space  $A_{\mu}^1$.

\BET \label{th1.1}
There are sequences  $0 < s_1 < s_2 < \ldots < R$ and $0=m_0 < m_1 < m_2 < \ldots $,
non-negative numbers $d_n$, $t_{n,k}$ (with $n \in \bbN$ and $[m_{n-1}]< k \leq
[m_{n+1}]$ ) and constants $c_1$, $c_2 >  0$
such that for all  $g(z) = \sum_{k=0}^{\infty}  \alpha_k z^k$  we have
\bea
c_1 \Vert g \Vert_1 \leq \sum_{n=0}^{\infty} M_1(T_ng,s_n)d_n \leq
c_2 \Vert g \Vert_1,   \label{1.8}
\eea
where
\bea
T_n g = \sum_{[m_{n-1}]+1}^{[m_{n+1}]}t_{n,k}  \alpha_k z^k.
\eea
\ENT

We will need the following consequence of this result.

\BEC \label{cor1.0}
Let $(n_j)_{j=1}^\infty$ be  an increasing sequence of indices such that
$n_{j+1} - n_j \geq 2 $ for all $j$ and let
$h_j(z) = \sum_{k= [m_{n_j}]+1}^{[m_{n_j+1}]} \alpha_ k z^k$
be a polynomial. We have
\bea
\Vert h \Vert_1 \leq \sum_{j=0}^{\infty} \Vert h_j \Vert_1 \leq
C \Vert h \Vert_1
\ \ \  
\mbox{for all}
\  
h = \sum_{j=1}^\infty h_j \in A_\mu^1 .  \label{1.8a}
\eea
\ENC

{\bf Proof.} Applying \eqref{1.8} to   $h_j$ yields that
$\Vert h_j \Vert_1 $ and
$$
M_1(T_{n_j} h_j,s_{n_j}) d_{n_j } + M_1(T_{n_j +1 } h_j,s_{n_j +1 })d_{n_j +1}
$$
are proportional quantities. Moreover, $T_nh = 0$, if $n$ is not equal to
$n_j$ or $n_j+1$ for  any $j$, and
\bea
T_{n_j} h = T_{n_j} h_j , \ \ T_{n_j +1 } h = T_{n_j +1 } h_j \ \ \ \mbox{for
all } j.
\eea
Hence, by another application of \eqref{1.8},
\beas
& & \Vert h \Vert_1 \leq  \sum_{j= 0}^\infty \Vert h_j\Vert_1
\leq C \sum_{j= 0}^\infty \big( M_1(T_{n_j} h_j,s_{n_j}) d_{n_j } + M_1(T_{n_j +1 } h_j,s_{n_j +1 })d_{n_j +1} \big)
\roweq
C \sum_{n= 0}^\infty M_1(T_n h,s_n) d_n \leq C' \Vert h \Vert_1 .
\ \hskip5cm \Box
\eeas

Let us make a remark concerning the numbers and constants in the
above results.

\BER
$1^\circ.$ Theorem \ref{th1.1} is a reformulation of Theorem 1.3. of
\cite{ahl}, where
the sequences $(s_n)_{n=1}^\infty$ and $(m_n)_{n=0}^\infty$ were
chosen, by using induction, such that, for all $n \in \bbN$,
\[
\int_0^{s_n} r^{m_n}d \mu = b \int_{s_n}^Rr^{m_n}d \mu \ \ \ \mbox{ and } \ \ \
\int_0^{s_n} r^{m_{n+1}}d \mu =  \int_{s_n}^Rr^{m_{n+1}}d \mu.
\]
where $b >5$ is some constant. Then, the numbers $d_n$ were set to be
\bea
d_n= \left(\int_0^{s_n}\left(\frac{ r}{s_n}\right)^{m_{n}}d \mu +
\int_{s_n}^R\left(\frac{r}{s_n}\right)^{m_{n+1}}d \mu \right) .
\eea
As proven in Section 5 of \cite{ahl}, it is always possible to find these
sequences, although calculating them exactly for given concrete weights
seems in general to be difficult.

$2^\circ$. If $R=1$ and  $d \mu = r v(r) dr d\theta$ with
$v(r) = \exp\big( - \alpha (1-r^\ell)^{-\beta}\big) $ for some constants
$\alpha, \beta, \ell > 0$, then the numbers $m_n$ and $s_n$, $n \in \bbN$,
($m_0 = 0$), were calculated by a different method than in $1^\circ$
in Propositions 3.1 and 3.3.$(ii)$  of \cite{BLT3}:
\bea
m_n = \ell \beta^2 \Big( \frac{\beta}{\alpha} \Big)^{1/\beta}
n^{2 + 2/\beta} - \ell \beta^2 n^2 \ \ \mbox{and} \ \
s_n = \Big( 1 -  \Big( \frac\alpha{\beta} \Big)^{1/\beta}
n^{- 2/\beta} \Big)^{1/ \ell}  .
\eea

$3^\circ$. In the citations mentioned in  $1^\circ$ and $2^\circ$, the numbers
$t_{n,k}$ were chosen as the coefficients of certain de la Valle\'e Poussin operators, more precisely,
\begin{equation}
t_{n,k}=
\left\{
\begin{array}{ll}
{ \displaystyle \frac{k - [m_n] }{[m_n]- [m_{n-1}]} } , &  \ \mbox{if} \
m_{n-1} < |k| \leq m_n , \\
 & \\
{ \displaystyle \frac{[m_{n+1}] - k }{[m_{n+1}]- [m_n]} },
&  \ \mbox{if}  \
m_{n} < |k|  \leq m_{n+1}  .
\end{array}
\right.
\end{equation}
\ENR

Let us next state our result on the characterization of solid
$A_\mu^1$ spaces.

\BET \label{th1.2}
 The following are equivalent:

(i) $A_{\mu}^1$ is solid,

(ii) $s(A_{\mu}^1) = A_{\mu}^1$,

(iii) the monomials $(z^n)_{n=0}^{\infty}$ are an unconditional basis of $A_{\mu}^1$,

(iv) the normalized monomials $(z^n/ \Vert z^n \Vert_1)_{n=0}^{\infty}$ are equivalent to the
unit vector basis of $\ell^1$,

(v)  $\sup\limits_{n \in \bbN}(m_{n+1}-m_n) < \infty$ for the numbers $m_n$ in Theorem \ref{th1.1}.
\ENT

In the following we retain the numbers $m_n$, $s_n$ of
Theorem \ref{th1.1} and consider the Dirichlet projections $P_n$.

\BEL
\label{lem1.4}
Assume that $\limsup_{n \rightarrow \infty}(m_{n+1}-m_n)= \infty$.
Then, for every $N> 0$ there exist  an arbitrarily large $n \in \bbN$, an index
$M < m_{n+1}$  and a
polynomial $f(z )= \sum_{k=[m_n]+1}^{[m_{n+1}]} \alpha_k z^k$
with $ \Vert f \Vert_1 \leq 1$ but $ \Vert P_{M}f \Vert_1
= \Vert ( P_{M} - P_{[m_n]} ) f \Vert_1 \geq N$.
\ENL

{\bf Proof.} Due to the unboundedness of the operator norms of $P_n$ on
$H^1$, see Section \ref{sec0}, we find an index $K$ and a polynomial $g(z)=
\sum_{j=0}^L \beta_j z^j$ with $M_1(g,1) = 1$ but $M_1(P_K g,1) > N$. By
assumption we find $n \in \bbN$, as large as we wish, such that
 $m_{n+1}-m_n > L+1$.  Then put
\[
f(z) = \sum_{k=[m_n]+1}^{[m_n]+L+1}\beta_{k-[m_n]-1}
\frac{1}{s_n^k}z^k.
\]
We obtain
\[
M_1(f,s_n) = M_1(g,1)=1 \ \ \ \mbox{ and }
\ \ \ M_1(P_{m_n+K}f, s_n) = M_1(P_Kg,1) > N.
\]
Put $M= K+ [m_n]+1$ and use Theorem \ref{th1.1} to complete the proof of the lemma. We have $P_M f = (P_{M} - P_{[m_n]}) f$ just by the choice of $f$.
\ \ $\Box$

\bigskip

{\bf Proof of Theorem \ref{th1.2}.} $(i) \Leftrightarrow (ii)$: follows from the definition. \\
$(iv) \Rightarrow (iii) \Rightarrow (ii)$: these are obvious. \\
$(ii) \Rightarrow (v)$: 	Assume that $\limsup_{n \rightarrow \infty}
(m_{n+1}-m_n)= \infty$. For every $j \in \bbN$ we find, by Lemma \ref{lem1.4},
a polynomial $f_j \in $ span$\{z^{[m_{n_j}]+1}, \ldots, z^{[m_{n_{j+1}]}}\}$
for some $m_{n_j}$ with
\bea
 \Vert f_j \Vert_1=2^{-j}  \ \ \mbox{and
$ \big\Vert  P_{k_j} f_j \big\Vert_1 \geq 1$
for some $k_j \in ( m_{n_j}, m_{n_j +1 } )$.  }   \label{1.8b}
\eea
We may assume that $n_{j+1}- n_j \geq 2$. Put $f= \sum_j f_j$
and $g = \sum_j P_{k_j}f_j = \sum_j ( P_{k_j} - P_{[m_{n_j}]} )f_j
$. Then, $f \in A_{\mu}^1$ but in view of \eqref{1.8b}, \eqref{1.8a} we have
$ g \not\in A_{\mu}^1$. Hence $f \not\in s(A_{\mu}^1)$.
\\
$(v) \Rightarrow (iv):$ Let $g(z) = \sum_{k=[m_{n-1}]+1}^{[m_{n+1}]} \alpha_kz^k$. By $(v)$
we obtain a constant independent of $n$, $r$ and $g$
with
\[ M_1(g,r) \leq \sum_{k=[m_{n-1}]+1}^{[m_{n+1}]} |\alpha_k| r^k \leq c M_1(g,r).\]
Then, \eqref{1.8} 
yields numbers $\delta_k = t_{n,k}s_n^k d_n$  such that for all functions
$f(z) =
\sum_{k=0}^{\infty} \alpha_k z^k \in A_\mu^1$ we have, with universal constants $c_1, c_2$,
\[
c_1 \Vert  f \Vert_1 \leq \sum_{k=0}^{\infty} \delta_k |\alpha_k| \leq c_2  \Vert f \Vert_1. \]
 This proves  $(iv)$.
$\Box$

\BEC
\label{cor1.5}
If $R < \infty$ then $A_{\mu}^1$ is never solid.
\ENC

{\bf Proof.} It follows from Proposition 2.1. of
\cite{ahl} that in this case we always have
$\limsup_{n \rightarrow \infty}(m_{n+1}-m_n)= \infty$. $\Box$
\bigskip

\BEE \label{ex1}  There are indeed examples where
$A_{\mu}^1$ is solid. Let $R= \infty$ and
$d\mu(r) = \exp(-\log^2(r))dr$. It was shown in \cite{ahl}, Example 2a) that here $ \sup_n (m_{n+1}-m_n) < \infty$.
\ENE

\BET
\label{th1.6}
Let $m_n$, $s_n$ and $d_n$ be the numbers of Theorem \ref{th1.1}. The solid core of $A_{\mu}^1$ equals
\bea s(A_{\mu}^1) = &\bigg\{  & g: R \cdot \mathbb{D} \rightarrow \mathbb{C} : g(z) = \sum_{k=0}^{\infty} \hat{g}(k)z^k  \nonumber
\\
& & \mbox{ with }
\sum_{n=1}^{\infty} d_n \Big( \sum_{k=[m_n]+1}^{[m_{n+1}]} |\hat{g}(k)|^2 s_n^{2k} \Big)^{1/2} < \infty \bigg\}.  \label{1.9c}
\eea
\ENT

{\bf Proof.} For a holomorphic function
$g(z) = \sum_{k=0}^{\infty} \hat{g}(k) z^k$
we write
\[ g_n(z) = \sum_{k=[m_n]+1}^{[m_{n+1}]}
\hat{g}(k)z^k \ \ \mbox{ and } \ \ g_I(z) = \sum_{n=0}^{\infty}g_{2n}(z), \ \ g_{II}(z) =\sum_{n=0}^{\infty} g_{2n+1}(z). \]
Let us denote by $V$ the function space on the right-hand side of \eqref{1.9c}.
Moreover, for all $n$, let $\Delta_n = \{ +1, -1 \}^{[m_{n+1}] -[m_n]}$, and
for $\Theta_n = (\theta_{[m_{n}]+1}, \ldots, \theta_{[m_{n+1}]})\in \Delta_n$ put
\[g_{\Theta_n}(z) = \sum_{k=[m_n]+1}^{[m_{n+1}]} \theta_k \hat{g}(k) z^k.  \]
At first assume that $g \in V$. Then $g_I, g_{II} \in V$. Let $f$ be
holomorphic with $|\hat{f}(k)| \leq |\hat{g_I}(k)|$ for all $k$. By
\eqref{1.8a} and Theorem \ref{th1.1}
\begin{eqnarray*}
  \Vert f \Vert_1 & \leq & \sum_{n=0}^{\infty}  \Vert f_{2n} \Vert_1 \leq
c \sum_{n=0}^{\infty} d_{2n} M_1(f_{2n}, s_{2n})
\\
& \leq &  c \sum_{n=0}^{\infty} d_{2n} M_2(f_{2n}, s_{2n})
\leq c \sum_{n=0}^{\infty} d_{2n} M_2(g_{2n}, s_{2n})
< \infty
\end{eqnarray*}
where $c >0$ is a universal constant and we also used the definition of the space $V$ in the last step. Hence $f\in A_{\mu}^1$, in particular
$g_I \in A_{\mu}^1$. We conclude $g_I \in s(A_{\mu}^1)$. The same proof shows that $g_{II}$ and hence  $g \in s(A_{\mu}^1)$.

Conversely, let $g \in s(A_{\mu}^1)$. Then
$g_I, g_{II} \in s(A_{\mu}^1)$.
Let ${\tilde{\Theta}}_n \in \Delta_n$ be such that
\[
a_1\left(\sum_{k=[m_n]+1}^{[m_{n+1}]}|\hat{g}(k)|^2\right)^{1/2} \leq \frac{1}
{2^{[m_{n+1}] -[m_n]}} \sum_{\Theta_n \in \Delta_n}M_1(g_{\Theta_n}, s_n)
\leq M_1(g_{{\tilde{\Theta}}_n} ,s_n).
\]
Here we used the Khintchine inequality (see \cite{Zy},  Ch.\,V, Thm.\,8.4)
with the Khintchine constant $a_1$. Put $h_I = \sum_{n=0}^{\infty} g_{
{\tilde{\Theta}}_{2n}}$. Then we obtain $|\hat{h_I}(k)| = |\hat{g_I}(k)|$  for all $k$.
Hence $h_I \in A_{\mu}^1$. The   choice of ${\tilde{\Theta}}_n$ and Theorem
\ref{th1.1} applied to $h_I$ yield
\beas
& & \sum_{n=0}^{\infty} d_{2n} \Big( \sum_{k=[m_{2n}]+1}^{[m_{2n+1}]} |
\hat{g}(k)|^2 s_{2n}^{2k} \Big)^{1/2} =
\sum_{n=0}^{\infty} d_{2n}
\Big( \sum_{k=[m_{2n}]+1}^{[m_{2n+1}]} |\hat{h_I}(k)|^2 s_{2n}^{2k}
\Big)^{1/2}
\rowleq
\frac{ 1}{a_1}\sum_{n=0}^{\infty} d_{2n}  M_1(g_{{\tilde{\Theta}}_{2n}} ,s_{2n}) \leq \frac{c_2}{a_1}  \Vert h_I \Vert_1 < \infty.
\eeas
Here, $c_2$ is the constant of Theorem \ref{th1.1}. We conclude $g_I \in V$,
and similarly we see that $g_{II} \in V$. Hence $g \in V$, which  implies $V =
s(A_{\mu}^1)$. \ \ $\Box$

\section{On solid hulls}
\label{sec2}
In this section we assume $R=1$. We start by the remark that in addition to the
definition of a solid hull as in Section \ref{sec0}, there exist two other a
priori different definitions in the literature: in \cite{AS}, the solid hull
$S_{{\rm vect}}(X)$ of a space $X$ of analytic functions on $\bbD$ is defined as
the intersection of all solid {\it  vector spaces} of analytic functions on
$\bbD$. Obviously, $ S(X)$ is a vector space if and only if for every $f,g \in X$
there  is $h \in X $ such that the Taylor coefficients satisfy
$|\hat f(k)|+|\hat g(k)| \leq |\hat h(k)|$ for all $k$.

One more variant appears in the theory of so called BK-spaces.
By definition, a BK-space is a vector space of complex sequences $f=
(f_k)_{k=0}^\infty$  endowed with a norm which makes it into a Banach space,
such that the coordinate functionals become bounded operators.
In the theory of BK-spaces, see \cite{Bu2}, the solid hull $S_{BK}(X)$ of a
BK-space $X$ is  defined as the intersection of all solid $BK$-spaces
containing $X$. By using Taylor coefficients
we consider Banach spaces of analytic functions on $\bbD$ as BK-spaces,
and, in particular, we will characterize in the sequel the solid hull $S_{BK}(A_\mu^1)$ although we
will  avoid using the terminology of BK-spaces, except for the proof of
Proposition \ref{prop2.9}. It is quite easy to see that
\bea
S(X) \subset S_{{\rm vect}}(X) \subset S_{BK} (X)
\eea
for a BK-space $X$ as above. All results on solid hulls $S(X)$ in the
literature, which are known to the  authors, happen to be vector spaces
which can be endowed with norms making them into solid BK-spaces. Thus, in all of
these cases one actually has $S(X) = S_{BK}(X)$.

Our aim is to use the known  duality relations between weighted
$A^1$ and $H^\infty$-spaces and existing results of the solid core  of
$H_v^\infty$ in order to find the solid hull $S_{BK} ( A_\mu^1) $.
We focus on the case the measure $\mu$ is the weighted Lebesgue measure $ v dA$
with a radial weight $v$ making the Bergman space into a "large" one:
the admissible weights include the exponentially decreasing weights,
see Example \ref{ex3.3}, below.

We start by some general considerations.

Given a sequence $\theta = (\theta_k)_{k=0}^\infty$ with $|\theta_k|\leq 1$ for
all $k$, we denote by $M_\theta$ the operator
$M_\theta  \sum_{k=0}^\infty \hat f(k)  z^k =
\sum_{k=0}^\infty \theta_ k \hat f(k) z^k$. We will need to consider analytic
function  spaces on $\bbD$  such that the norm of the space satisfies
\bea
\Vert M_\theta f \Vert \leq \Vert f \Vert  \label{4.4}
\eea
for all $f = \sum_{k=0}^\infty \hat f(k) z^k\in X $ and all sequences $\theta = (\theta_k)_{k=0}^\infty$
with $|\theta_k|\leq 1$ for all $k$.

The following result is essentially known.

\BEP \label{prop2.9}
If $( X, \Vert \cdot \Vert_X ) $ is a Banach space of analytic functions on the
unit disc $\bbD$ such that all coordinate functionals $ f \mapsto \hat f(k)$
are bounded operators, then its solid  hull  $S_{BK}(X)$ can be endowed with a norm
$\Vert \cdot \Vert_S$ such that

\noindent $(i)$ the embedding $X \hookrightarrow S_{BK}(X)$ is continuous,

\noindent $(ii)$ the norm $\Vert \cdot \Vert_S$ satisfies \eqref{4.4},

\noindent $(iii)$ if $p : S_{BK}(X) \to \bbR_0^+$ is any norm with \eqref{4.4}
such that
$p(f) \leq \Vert f \Vert_X$ for all $f \in X$, then $p(f)  \leq C \Vert f
\Vert_S$  for a constant $C > 0$ and all $f \in S_{BK}(X)$,

\noindent $(iv)$ the normed space $\big( S_{BK}(X), \Vert \cdot \Vert_S \big)$
is complete, and

\noindent $(v)$ if the subspace of polynomials $\cP$ is dense in
$X$, then it is dense in $\big( S_{BK}(X), \Vert \cdot \Vert_S \big)$,
too.
\ENP

Proof.  Let us explain how the claims follow from the theory of BK-spaces, see
\cite{Bu1}, \cite{Bu2}. For the sake of the simplicity of notation, let us
consider $X$  as a BK-sequence space in the following, which we can do by
assumption. We denote by $y \cdot f$ the  coordinatewise
product of two complex sequences $y$ and $f$. The space $
\ell^\infty \widehat \otimes X $ is defined in  \cite{Bu1} to
consist of sequences $g= (g_k)_{k=0}^\infty$ having a coordinatewise
convergent representation
\bea
g = \sum_{j=1}^\infty y^{(j)} \cdot f^{(j)} \ \ \ \mbox{with
 $y^{(j)} = \big( y^{(j)}_k \big)_{k=0}^\infty \in \ell^\infty$,
  $f^{(j)} = \big( f^{(j)}_k \big)_{k=0}^\infty \in X \ \forall\,j$}
\label{4.5a}
\eea
such that
\bea
\sum_{j=1}^\infty  \Vert y^{(j)} \Vert_{\ell^\infty} \Vert f^{(j)} \Vert_X
< \infty . \label{4.5b}
\eea
The  norm $\Vert \cdot \Vert_S$ of $g \in \ell^\infty \widehat \otimes X$ is defined by taking the
infimum of the quantity \eqref{4.5b} over all possible representations
\eqref{4.5a} of $g$. Theorem 3 of \cite{Bu1} yields that the resulting space
is complete, and Theorem 8 of \cite{Bu2}  says that $\ell^\infty \widehat
\otimes X$ equals the solid hull $S_{BK}(X)$. The completeness of the
space is included in the same reference, hence, property $(iv)$ holds.

If $f \in X$, then we have ${\rm e} \cdot f = f$, where ${\rm e} =
(1,1,1,\dots )  \in \ell^\infty$, and in view of the above definition of the
norm $\Vert \cdot \Vert_S$, this implies that $\Vert f \Vert_S \leq \Vert f \Vert_X$ for all $f \in X$ so
that  the embedding of $X$ into $\big( S_{BK}(X), \Vert \cdot \Vert_S
\big) $ is
continuous.

Also,
if $g \in S_{BK}(X) $ has a representation \eqref{4.5a} and $\theta$ is
given as in \eqref{4.4}, then $M_\theta g$ has a coordinatewise convergent
representation
\bea
M_\theta g = \sum_{j=1}^\infty ( M_\theta y^{(j)} ) \cdot f^{(j)} ,
\eea
and property $(ii)$ follows from the definition of $\Vert \cdot \Vert_{BK}$.

In the proof of Theorem 3 of \cite{Bu1} it is shown if $p$ is  the norm of
any BK-space containing $S_{BK}(X)$, then there exists $C> 0$ such that
\beas
p(y \cdot f) \leq C \Vert y \Vert_\infty \, \Vert f \Vert_X
\eeas
for all $y \in \ell^\infty$, $ f \in X$. This implies
\beas
p\Big( \sum_{j=1}^\infty y^{(j)} \cdot f^{(j)} \Big) \leq
C  \sum_{j=1}^\infty \Vert y^{(j)}  \Vert_\infty \, \Vert f^{(j)} \Vert_X
\ \ \mbox{for all} \ g =  \sum_{j=1}^\infty y^{(j)} \cdot f^{(j)}
\in \ell^\infty \widehat \otimes X,
\eeas
and property  $(iii)$ follows from the definition of $\Vert \cdot \Vert_S$.
Finally, as for property  $(v)$, it follows from Theorem 2 of \cite{Bu1}
that finite linear combinations of functions $y \cdot f$, $y \in \ell^\infty$,
$f \in X$, form a dense subspace of $\ell^\infty \widehat \otimes X = S(X)$.
If $y$, $f$ and $\varepsilon > 0$ are given, we use the assumption
in $(v)$ to find a polynomial $h$ such that
$\Vert f- h\Vert_X < \varepsilon / (1 + \Vert y \Vert_\infty)$. Then,
$y \cdot h$ is a polynomial, which  satisfies
\beas
p( y \cdot f - y \cdot h) = p( y \cdot ( f-h))
\leq \Vert y \Vert_\infty \Vert f-h \Vert_X \leq \varepsilon .
\eeas
Property $(v)$ follows from these arguments. \ \ $\Box$

\BEL
\label{lem2.0}
Let $X$  be a Banach space of analytic functions on the unit disc $\bbD$
such that the subspace $\cP$ of polynomials is dense in $X$, and
let $w$ be a radial weight function on $\bbD$.
Let  $Y$ be the  space of all analytic functions  $g$ on the
disc such that
\bea
\sup\limits_{f \in B_X} |\langle f , g \rangle| < \infty, \ \ \
\mbox{where} \ \ \ \
\langle f , g \rangle = \int\limits_\bbD f \overline{g} w dA   \label{2.a}
\eea
and $B_X$ denotes the unit ball of $X$.
If $X$ is solid and there exists a constant $C > 0$ such that
\bea
\Vert M_\theta f \Vert_X \leq C \Vert f \Vert_X
\label{2.b}
\eea
for all numerical sequences $\theta = (\theta_k)_{k=0}^\infty$ with $|\theta_k| \leq 1$,
then $Y$ is solid, too.
\ENL

We point out  given a  Banach space $X$ as in the assumption, it is not in general
known whether its dual space has a representation as a space of analytic functions
with dual norm coming from \eqref{2.a}.

\bigskip

Proof. If $g = \sum_{k=0}^\infty \hat g(k)  z^k \in Y$ and  $\theta$
is as above, then for
$ 
M_ \theta g 
$ 
we have by \eqref{2.b}
\bea
& & 
\sup\limits_{f \in B_X}
|\langle f , M_\theta g  \rangle |
= \sup\limits_{f \in B_X} \sum_{k=0}^\infty \bar \theta_k
\hat f(k)    \overline{ \hat g(k)} \int\limits_0^1 r^{2k+1} w(r)dr
\roweq
\sup\limits_{f \in B_X}
|\langle M_{\bar \theta} f ,  g  \rangle |
\leq
\sup\limits_{\stackrel{ \scriptstyle f \in X }{\Vert f \Vert_X \leq C } }
|\langle f , g  \rangle | < \infty . 
\label{2.e}
\eea
Thus, $M_\theta g \in Y$. \ \ $\Box$

\bigskip

We next recall an elementary fact concerning Banach sequence
spaces. Assume that the sequences $(\beta_k)_{k=0}^\infty$ and
$(\gamma_k)_{k=0  }^\infty$ of positive numbers are given
and  $\alpha_k = \gamma_k \beta_k^{-1}$ for all $k$. Let also
$(\mu_n)_{n=0}^\infty$  be an increasing, unbounded sequence of non-negative
numbers; denote $\mu_{-1} = -1$ and let
\bea
A &=& \big\{ a = (a_k)_{k=0}^\infty \, : \,
\Vert a \Vert_A = \sum_{n \in \bbN}  \max_{\mu_{n-1} < k \leq  \mu_n }
\alpha_k |a_k|  < \infty \big\},
\label{2.24c}\\
B &=& \big\{ b = (b_k)_{k=0}^\infty \, : \,
\Vert b \Vert_B = \sup\limits_{n \in \bbN}  \sum_{\mu_{n-1} < k \leq \mu_n }
 \beta_k |b_k| < \infty \big\} .  \label{2.24a}
\eea
Then,  $B$ is the dual of $A$ with respect to the dual pairing
\bea
\langle a, b \rangle = \sum_{k=0}^\infty  \gamma_k b_k \overline{a_k},
\ \ \ \mbox{where} \ a=(a_k)_{k=0}^\infty \in A ,
\ b=(b_k)_{k=0}^\infty \in B .
\label{2.24b}
\eea

From now on we consider radial weights  $v : \bbD \to \bbR^+$
satisfying the following two assumptions.

\medskip

\noindent $(I)$ We have
\bea
v(z) = \exp ( - \varphi(z) )  ,
\eea
where $\varphi$ belongs to the class $\cW_0$ of \cite{HLS}.

\medskip

\noindent We will not need a detailed definition of $\cW_0$,
but recall that $\varphi \in \cW_0$, if it is
a twice continuously differentiable real valued function  with
$\Delta \varphi > 0$ on $\bbD$ and there exists a  function
$\rho : \bbD  \to  \bbR$ and a constant $C >0$ such that
\bea
\frac{1}{C}\rho (z) \leq \frac{1}{\sqrt{\Delta \varphi(z)}} \leq C \rho(z)
\ \ \forall \, z \in \bbD ;
\eea
the function $\rho $ must also satisfy the H\"older-property
\bea
\sup\limits_{z,w \in \bbD, z\not=w} \frac{|\rho(z) - \rho(w)|}{|z-w|}
< \infty
\eea
as well as the Lipschitz-property
\bea
\forall \, \varepsilon > 0 \, \exists \,  \mbox{compact}\,E \subset \bbD:
\  |\rho(z) - \rho (w) | \leq \varepsilon |z-w| \ \forall \, z,w
\in \bbD \setminus E.
\eea
For more details, see \cite{HLS}. Note that the considerations in
\cite{HLS} are not restricted to radial weights, contrary to our
situation.

According to \cite{HLS}, Theorem 4.3., if the weight $v$ satisfies condition
$(I)$, then  the space $H_v^\infty$   is the dual of   $A^1_v$ with respect to the dual
pairing
\bea
\langle f , g \rangle = \int\limits_\bbD f \overline g v^2 dA  .\label{2.4}
\eea

The second requirement is the following:
\medskip

\noindent $(II)$ The weight $v$ satisfies the condition $(b)$ of
\cite{BT1}, \cite{BLT1}.

\medskip

\noindent Recall that the weight $v$ satisfies the condition
$(b)$ if there exist numbers $b> 2$, $K > b$ and $ 0 < \mu_1 < \mu_2 < \ldots$ with
$\lim_{n \rightarrow \infty} \mu_n = \infty$ such that
\bea
b \leq \left( \frac{r_{\mu_n}}{r_{\mu_{n+1}}} \right)^{\mu_n}
\frac{v(r_{\mu_n})}{v(r_{\mu_{n+1}})}, \left( \frac{r_{\mu_{n+1}}}{r_{\mu_{n}}} \right)^{\mu_{n+1}}
\frac{v(r_{\mu_{n+1}})}{v(r_{\mu_{n}})} \leq K ,     \label{bbb}
\eea
where $r_m \in ]0,1[ $ denotes the global maximum point of the function
$r^m v(r)$ for any $m > 0$. Theorem 2.4 of  \cite{BLT1} states that the solid
core of the space  $H_v^\infty$ equals
\bea
& &  s(H_v^\infty)
= 
\Big\{ (b_k)_{k=0}^\infty \, : \, \Vert b \Vert_{v,s} = \sup\limits_{n \in \bbN} v(r_{\mu_n})
\sum_{\mu_n < k \leq \mu_{n+1}} |b_k| \sigma_k <
\infty \Big\},   \label{2.23}
\eea
where we denote $\sigma_k = r_{\mu_n}^k$.  Let us define for every  $k  \in \bbN_0$ the number
\bea
S_k = \frac{\int\limits_0^1 r^{2k+1} v(r)^2dr}{v(r_{\mu_n}) \sigma_k} \  , \label{2.24}
\eea
where $n$ is the unique number such that $\mu_n < k \leq \mu_{n+1}$.

\BEE
\label{ex3.3}
According to \cite{BLT1}, all weights $v(r) = \exp \big( - \alpha /
(1-r^2)^\beta \big)$ with  $\alpha, \beta>0$, satisfy condition $(b)$, and it
is easy to see that they  also satisfy  assumption $(I)$.
\ENE

\BET
\label{cor2.3}
Let the weight $v$ satisfy the assumptions $(I)$ and $(II)$.
Then, we have
\bea
& &  S_{BK}(A_{\mu}^1)
= 
\Big\{ b = (b_k)_{k=0}^\infty \, : \, \Vert b \Vert_{\mu,S} =
\sum_{n=0}^{\infty}
\sup_{\mu_n < k \leq \mu_{n+1}} |b_k| S_k < \infty \Big\},   \label{2.25}
\eea
and the norm $\Vert \cdot \Vert_S$ given by Proposition \ref{prop2.9} is equivalent
with  $\Vert \cdot \Vert_{\mu, S}$.
\ENT

Proof. Let the solid hull $S_{BK}(A_\mu^1) $ be endowed with  the norm
$\Vert \cdot \Vert_S$  of Proposition \ref{prop2.9}, and let us denote
the Banach space on the right-hand side of \eqref{2.25} by $Z$.

We note that by
the duality relations explained above (see \eqref{2.23} for the
definition of $\Vert \cdot \Vert_{v,s}$), we have for all $f \in A_\mu^1$
\bea
\Vert f \Vert_1 = \sup\limits_{\stackrel{\scriptstyle g \in H_v^\infty}{
\scriptstyle \Vert g \Vert_{H_v^\infty} \leq 1}} |\langle f , g \rangle|
\ \ \mbox{and} \ \
\Vert f \Vert_{\mu,S} = \sup\limits_{\stackrel{\scriptstyle g \in
s(H_v^\infty)}{
\scriptstyle \Vert g \Vert_{v,s} \leq 1}} |\langle f , g \rangle| .
\eea
It is proved in \cite{BLT1}, equation (2.4) and the very end of the proof of
Theorem 2.4, that $\Vert g \Vert_{H_v^\infty} \leq C \Vert g \Vert_{v,s}$
for $g \in s(H_v^\infty)$. Therefore $\Vert f \Vert_1 \geq C \Vert f \Vert_{\mu,S}$ for all
$f \in A_\mu^1$. This implies in particular that $A_\mu^1 \subset Z$.
Clearly, $Z$ is a solid Banach space
and the coordinate functionals are continuous, thus
it contains the space $S_{BK}(A_\mu^1)$. Moreover, we obtain
$\Vert f \Vert_S \geq C \Vert f \Vert_{\mu,S}$
for $f \in S_{BK}(A_\mu^1)$ from Proposition \ref{prop2.9}.$(iii)$.

We show that the
norms $\Vert \cdot \Vert_{\mu,S}$ and $ \Vert \cdot \Vert_S$ are equivalent
in $S_{BK}(A_\mu^1)$. To do this, we  prove that  $C \Vert f \Vert_{\mu,S}  \geq  \Vert f \Vert_S$.
Note  that the space \eqref{2.23} is the dual space of
$Z$ 
in the dual pairing \eqref{2.4}. Indeed, if
$f=\sum_k \hat f(k)  z^k$
and $g  = \sum_k \hat g(k) z^k $ 
are polynomials, then, by a direct calculation,
\bea
\langle f, g \rangle = \sum_{k=0}^\infty \hat f(k) \overline{\hat g(k)}
\int\limits_0^1 r^{2k +1} v(r)^2 dr .  \label{2.39}
\eea
The result follows from \eqref{2.24c}--\eqref{2.24a},
in addition to the definitions \eqref{2.4}--\eqref{2.25}.

Suppose now by antithesis that $\Vert \cdot \Vert_{\mu,S} $ and
$\Vert \cdot \Vert_S$ are non-equivalent norms so that we can find
a sequence $(f_n)_{n=1}^\infty \subset S_{BK}(A_\mu^1)$ such that
\bea
\Vert f_n \Vert_{\mu,S}  \leq 2^{-n} \Vert f_n \Vert_S \ \ \mbox{and} \ \
\Vert f_n \Vert_S = 1 \ \forall \, n \in \bbN. \label{2.40}
\eea
By property $(v)$ in Proposition \ref{prop2.9} we can assume that $f_n$'s are
polynomials. We claim that it is possible to
find polynomials  $\tilde f_n$, $n \in \bbN$, with property \eqref{2.40} such that they
have distinct degrees,  more precisely
\bea
\tilde f_n(z) = \sum_{k=K_n}^{K_{n+1} -1 } \hat f (n,k) z^k , \ \ n \in \bbN,
\eea
for some unbounded sequence $0 = K_0 < K_1 < \ldots$ and some $\hat f (n,k)
\in \bbC$. Assume that $N \in \bbN$ and that such polynomials $\tilde f_n$ have
been found for
$n \leq N$, and let $M \in \bbN$ be the highest degree of these polynomials. Since
$\cP_M$ (the $M+1$-dimensional space of polynomials of degree at most $M$)
is  finite dimensional, all norms are equivalent there and we thus
find a constant $K = K(M) > 0$ such that
\bea
\Vert f \Vert_S \leq K  \Vert f \Vert_{\mu,S}  \label{2.42}
\eea
for all $f \in \cP_M$. We pick up the polynomial $f_L$  as in \eqref{2.40}
with $L = M+K$ and write $ f_1 = P_M f_L $, $f_2 = f_L- f_1$,
where $P_M$ is the $M$th Dirichlet projection from $S_{BK}(A_\mu^1)$ onto $\cP_M$, see
Section \ref{sec0}. Then, we have
$\Vert f_2\Vert_S \geq \frac12 \Vert f_L \Vert_S$,
since otherwise we get by  \eqref{2.42} and the triangle inequality
\beas
\Vert f_L \Vert_{\mu,S} \geq  \Vert f_1\Vert_{\mu,S}
\geq \frac{1}{K} \Vert f_1\Vert_S \geq \frac{1}{2K} \Vert f_L \Vert_S
 >  \frac{1}{2L}  \Vert f_L \Vert_S
\eeas
which contradicts with \eqref{2.40}. Now we get
\bea
\Vert f_2 \Vert_{\mu,S} \leq \Vert f_L \Vert_{\mu,S}
\leq 2^{-L}  \Vert f_L \Vert_S \leq 2^{-L + 1 } \Vert f_2 \Vert_S.
\eea
Taking $f_2 \Vert f_2\Vert_S^{-1}$ for $\tilde f_{N+1}$, the claim is proved.

Finally, for every $n$ we set
\bea
T_n := \Big(P^{(n)}   \big( S_{BK}(A_\mu^1) \big), \Vert \cdot \Vert_S \Big) \ \ \mbox{with} \ \
P^{(n)} = P_{K_{n+1}-1} - P_{K_n}
\eea
and then, using the Hahn-Banach theorem, pick up a polynomial
$$
g_n = \sum_{k=K_n}^{K_{n+1} -1 } \hat g (n,k) z^k
$$
which  defines a bounded functional on $(T_n, \Vert \cdot \Vert_S)$ with
respect to the dual pairing \eqref{2.39}, such that
\bea
\langle \tilde f_n , g_n \rangle = 1 , \ \ \Vert g_n \Vert_{n,*}
:= \sup\limits_{\stackrel{\scriptstyle f \in T_n}{
\scriptstyle \Vert f \Vert_S \leq 1}} |\langle f , g_n \rangle|= 1   \label{2.44}
\eea
Then, we observe that $g_n$ extends via \eqref{2.39}
to a functional on $S_{BK}(A_\mu^1)=: S$ such that
\bea
\sup\limits_{\stackrel{\scriptstyle f \in S}{
\scriptstyle \Vert f \Vert_S \leq 1}} |\langle f , g_n \rangle|=
\sup\limits_{\stackrel{\scriptstyle f \in S}{
\scriptstyle \Vert f \Vert_S \leq 1}} |\langle P^{(n)} f , g_n \rangle|=
\sup\limits_{\stackrel{\scriptstyle f \in T_n}{
\scriptstyle \Vert f \Vert_S \leq 1}} |\langle f , g_n \rangle| = 1,  \label{2.46}
\eea
since the norm $\Vert \cdot \Vert_S$ of $S_{BK}(A_\mu^1)$ satisfies
$(ii)$ of Proposition \ref{prop2.9} and  thus
$\Vert P^{(n)}  f \Vert_S \leq \Vert f \Vert_S$ for all $f \in S_{BK}(A_\mu^1)$.
Consequently,
\bea
g = \sum_{n \in \bbN} \frac{1}{n^2} g_n   \label{2.47}
\eea
is an analytic function which also is a bounded functional on
$\big( S_{BK}(A_\mu^1), \Vert \cdot \Vert_S \big) $ in
the dual pairing \eqref{2.39}. However, $g$ is not a bounded functional on $Z$,
since
\bea
\langle 2^n \tilde f_n , g \rangle = \frac{2^n}{n^2} \langle \tilde f_n , g_n \rangle
= \frac{2^n}{n^2}
\eea
and by \eqref{2.40}, the $Z$-norm $\Vert 2^n \tilde f_n \Vert_{\mu,S}$  is still at most 1.

The space $Y$ of all analytic functions on $\bbD$, which also are bounded
functionals on $\big(S_{BK}(A_\mu^1), \Vert \cdot \Vert_S \big)$ in the dual pairing
\eqref{2.39},  equals the space $Y$ in Lemma \ref{lem2.0}, when $X:= S_{BK}
(A_\mu^1)$.  Hence, $Y$ is solid. Due to the characterization of $H_v^\infty$ as
the dual of $A_v^1$, see  \eqref{2.4}, we also have $Y \subset H_v^\infty$.
On the other hand, we observed in the beginning of the proof that the solid core
$s(H_v^\infty)$, see \eqref{2.23}, equals the dual of $Z$ in the pairing
\eqref{2.39}. The properties of the function $g$, \eqref{2.47},  show that
$s(H_v^\infty) \subsetneq Y$, which contradicts the definition of a solid core.
We conclude that $C \Vert f \Vert_{\mu,S} \geq \Vert f \Vert_S$
for all $f \in S_{BK}(A_\mu^1)$.

We come to the conclusion that the norms $\Vert \cdot \Vert_S$ and
$\Vert \cdot \Vert_{\mu,S}$ are equivalent, hence, the spaces
$S_{BK}(A_\mu^1)$ and $Z$ coincide, since they both are complete.
\ \ $\Box$

\bigskip

Statement on Compliance with Ethical Standards:

\medskip

\noindent The research of Bonet was partially supported by the
project MCIN PID2020-119457GB-I00/AEI/10.13039/501100011033.

\noindent
There are no conflicts of interest as regards to this article.

\noindent
The research has been ethically conducted.

\noindent
The research does not involve human participants and/or animals.

\noindent
Data Availability Statements: n/a

\end{document}